\documentclass[10pt]{article}
\usepackage{latexsym}
\usepackage{amssymb,amstext,amssymb,amscd,euscript}
\usepackage{epsfig}
\usepackage{graphicx}
\usepackage{colortbl}
\usepackage{multicol}
\definecolor{lichtgrijs}{rgb}{0.9,0.9,0.9}
\definecolor{grijs}{rgb}{0.8,0.8,0.8}

\usepackage{theorem}
\theorembodyfont{\upshape}
\newtheorem{theorem}{Theorem}[section]
\newtheorem{definition}[theorem]{Definition}
\newtheorem{remark}[theorem]{Remark}

\newtheorem{lemma}[theorem]{Lemma}

\newtheorem{proposition}[theorem]{Proposition}
\newtheorem{corollary}[theorem]{Corollary}

\newcommand{\beq}{\begin{equation}}
\newcommand{\eeq}{\end{equation}}

\newcommand{\qed}    {\hspace*{\fill} $\Box$ \medskip}

\newcommand{\x}{\mathbf x}
\newcommand{\e}{\mathbf e}
\newcommand{\y}{\mathbf y}
\newcommand{\z}{\mathbf z}
\newcommand{\0}{\mathbf 0}
\newcommand{\T}{\ensuremath{\mathcal{T}}}
\newcommand{\V}{\ensuremath{\mathcal{V}}}
\newcommand{\W}{\ensuremath{\mathcal{W}}}
\newcommand{\C}{\ensuremath{\mathcal{C}}}
\newcommand{\X}{\ensuremath{\mathcal{X}}}
\newcommand{\B}{\ensuremath{\mathcal{B}}}
\newcommand{\mod}{{~\mathrm{mod}~}}

\title{Limiting shapes for deterministic centrally seeded growth models}
\author{Anne Fey-den Boer\footnote{Vrije Universiteit,
De Boelelaan 1081a,
1081 HV Amsterdam, The Netherlands, fey@eurandom.tue.nl} \and Frank Redig\footnote{Mathematisch Instituut Universiteit Leiden, Niels Bohrweg 1,
2333 CA Leiden,
The Netherlands, redig@math.leidenuniv.nl}}

\begin{document}
\maketitle

\begin{abstract}
We study the rotor router model and two deterministic sandpile models. For the rotor router model in $\mathbb{Z}^d$, Levine and Peres proved that the limiting shape of the growth cluster is a sphere. For the other two models, only bounds in dimension 2 are known. A unified approach for these models with a new parameter $h$ (the initial number of particles at each site), allows to prove a number of new limiting shape results in any dimension $d \geq 1$.

For the rotor router model, the limiting shape is a sphere for all values of $h$. For one of the sandpile models, and $h=2d-2$ (the maximal value), the limiting shape is a cube. For both sandpile models, the limiting shape is a sphere in the limit $h \to -\infty$. Finally, we prove that the rotor router shape contains a diamond.
\end{abstract}

\section{Introduction}

In a growth model, there is a dynamical rule by which vertices of a graph are added to an initial collection. The existing literature on growth models deals mainly with stochastic growth models. Among many stochastic growth models, there are for example the Eden model, the Richardson model, first and last passage percolation and diffusion limited aggregation. An introduction to stochastic growth models and limiting shape theorems can be found in \cite{durrett}.

An example related to the models in this paper, is the internal diffusion limited aggregation (IDLA) model. One adds particles one by one to the origin, letting each particle perform a random walk that stops when it hits an empty site. The growth cluster is then the collection of sites that contain a particle.

The IDLA model on $\mathbb{Z}^d$ has been studied by Lawler, Bramson and Griffeath. In their paper from 1992 \cite{lawler}, they prove that the limiting shape of the growth cluster for this model in any dimension is a sphere. Lawler \cite{lawler2} estimated the speed of convergence.

The three deterministic models to be discussed in this paper can be viewed as deterministic analogues of IDLA, and are all examples of a general height-arrow model studied in \cite{ham}. The closest analogue is the rotor router model, proposed by Propp (see \cite{kleber}). Once more, particles are added to the origin of the $d$-dimensional lattice $\mathbb{Z}^d$, but now they perform a deterministic walk as follows: at each site, there is an arrow present, pointing at one of the $2d$ neighboring sites. If a particle at this site finds it already occupied, then it takes a step in the direction of the arrow, and the arrow rotates to the next position.

In our other models, occupied sites hold the particles until they can be sent out in pairs, in the directions of a two-pointed arrow (double router model), or groups of $2d$, one in every direction (sandpile model). In the last model, exactly the same amount of particles is sent in every direction. In the other two models, the difference is at most one.

The last two models are deterministic versions of sandpile models. A stochastic sandpile model, where addition occurs at a random site, has originally been proposed to study self-organized criticality \cite{bak}. Soon after, Dhar \cite{dhar} introduced the abelian sandpile model, that has the symmetric toppling rule (``toppling'' of a site means that particles are sent out), which is now the most widely studied variant. Manna's sandpile model \cite{manna} is doubly stochastic: the addition site is random, and moreover, in a toppling two particles are sent to randomly chosen neighbors.

The rotor router model on $\mathbb{Z}^d$ has been studied by Peres and Levine \cite{levine,peres}. They have found that the limiting shape of this model is also a sphere, and give bounds for the rate of convergence. This result is at the basis of a number of new results in the present paper. Recently, a new paper of these authors appeared \cite{peres2}, further extending their results.

The deterministic abelian sandpile model has been studied on a finite square grid \cite{lubeck, markosova, wiesenfeld}, with emphasis on the dynamics, but hardly as a growth model. Le Borgne and Rossin \cite{rossin} found some bounds for the limiting shape for $d=2$.

From the above introductory description, it should be clear that all the deterministic models we introduced here are closely related.
In Section \ref{definitionssection}, we define the above models on $\mathbb{Z}^d$ in a common framework, to allow comparisons of the models. We introduce a parameter $h$ to parametrize the initial configuration. For all models, a (large) number $n$ of particles starts at the origin, and spreads out in a deterministic way through topplings. The rest of the grid initially contains a number $h$ of particles at each site, where $h$ can be negative. One can imagine a negative amount of particles as a hole that needs to be filled up; admitting negative particle numbers will be helpful in comparing the different models. Once every site is stable, (i.e., has a number of particles at most some maximal allowed number), the growth cluster is formed.

We present pictures, obtained by programming the models in Matlab, of growth clusters for finite $n$ and $d=2$. For all models, we obtain beautiful, self-similarly patterned shapes. The appeal of these pictures has been noted before, e.g., the sandpile pictures for $h= 0,-1$ can be found on the internet as ``sandpile mandala'', but the patterns are so far unexplained.
Similar patterns are found in the so-called sandpile identity configuration \cite{rossin,liu}.

Finally, we explain our main limiting shape results for each model. For the sandpile model, we obtain that the limiting shape is a cube for $h = 2d-2$, (Theorem \ref{cube}), and a sphere as $h \to -\infty$ (Theorem \ref{hminoneindig}). This last result also applies for the double router model. We remark that the sandpile model with $h \to -\infty$, strongly resembles the divisible sandpile introduced by Levine and Peres \cite{peres2}, of which the limiting shape is also a sphere.

For the rotor router model, we find that the limiting shape is a sphere for all $h$ (Theorem \ref{RRshape}).
Finally, we generalize the bounds of Le Borgne and Rossin for growth cluster of the sandpile model, to all $h$ and $d$. As a corollary, we find that the rotor router shape contains a diamond, and is contained in a cube.

The rest of the paper contains proofs of these results. In Section \ref{comparemodels}, we derive inequalities for the growth clusters. In Section \ref{shaperesults}, we prove the various limiting shape theorems.

\section{Model definitions and results}
\label{definitionssection}

Before introducing each of the three models, we present some general definitions that apply to each model. All models are defined on $\mathbb{Z}^d$.
We define a configuration $\eta = (H,T,D)$ as consisting of the following three components: the particle function $H:\mathbb{Z}^d \to \mathbb{Z}$, the toppling function $T:\mathbb{Z}^d \to \mathbb{Z}$ and the direction function $D:\mathbb{Z}^d \to \mathcal{D} = \{0,1, \ldots,2d-1\}$. We will use the $2d$ possible values in $\mathcal{D}$ to indicate the $2d$ unit vectors $\e_0 \ldots \e_{2d-1}$. The results will not depend on which value is assigned to which vector, as long as the assignment is fixed.

We call a configuration {\em allowed} if
\begin{enumerate}
\item For all $\x$, $T(\x) \geq 0$,
\item For all $\x$ with $T(\x)>0$, $H(\x)\geq 0$,
\item For all $\x$ with $T(\x)=0$, $H(\x) \geq h$,
\item For all $\x$, $D(\x) \in \mathcal{D}$.
\end{enumerate}
Here, $h \in \mathbb{Z}$, the ``background height'', is a model parameter.

Each model starts with initial configuration $\eta_0 = \eta_0(n,h)$, which is as follows: $H_0(\x) = n$ if $\x = \mathbf{0}$, and $H_0(\x) = h$ otherwise, $T_0(\x) = 0$ for all $\x$, and $D_0(\x) \in \mathcal{D}$ for all $\x$. Observe that $\eta_0$ is allowed.

We now define a toppling as follows:

\begin{definition}
A toppling of site $\x$ in configuration $\eta$ consists of the following operations:
\begin{itemize}
\item
$T(\x) \rightarrow T(\x) + 1$,
\item
$H(\x) \rightarrow H(\x) - c$, with $c \leq 2d$ some value specific for the model,
\item
$H(\x+\e_i) \rightarrow H(\x+\e_i) + 1$, with $i = (D(\x)+1) \mod 2d, \ldots, (D(\x)+c) \mod 2d$,
\item
$D(\x) \to (D(\x) + c)\mod 2d$.
\label{topplingdef}
\end{itemize}
\end{definition}

In words, in a toppling $c$ particles from site $\x$ move to $c$ different neighbors of $\x$, chosen according to the value of $D(\x)$ in cyclic order. We call a toppling of site $\x$ {\em legal} if after the toppling $H(\x) \geq 0$.
We call a site $\x$ stable if $H(\x) \leq h_{\max}$, with $h_{\max} = c-1$, so that a site can legally topple only if it is unstable.

We can now define stabilization of a configuration as performing legal topplings, such that a stable configuration is reached, that is, a configuration where all sites are stable. We call this configuration the final configuration $\eta_n = \eta_n(n,h)$, which will of course depend on the model. To ensure that the final configuration is reached in a finite number of topplings, we impose for each model $h<h_{\max}$. For each model, it is then known that $\eta_n$ does not depend on the order of the topplings. This is called the abelian property. The abelian property has been proved for centrally seeded growth models in general (\cite{diaconis}, Section 4), and for the sandpile model in particular \cite{meester}.

In the remainder of this paper, we will choose various orders of topplings. In Section \ref{RRshapesection}, we will even admit illegal topplings. In Section \ref{topplingcluster}, we will organize the topplings into discrete time steps, obtaining after time step $t$ a configuration $\eta^t$, consisting of $H^t$, $T^t$ and $D^t$.

We define two growth clusters that are formed during stabilization, denoting by $\x^\square$ the unit cube centered at $\x$ (i.e., $\x^\square = \{\y:\y = \z+\x, \z \in [-\frac 12, \frac 12]^d\}$):
\begin{definition}
The toppling cluster is the cluster of all sites that have toppled, that is,
$$
\T_n = \bigcup_{\x:T_n(\x) > 0} \x^\square,
$$
and the particle cluster is the cluster of all sites that have been visited by particles from the origin, that is,
$$
\V_n =  \T_n ~\cup ~\bigcup_{\x: H_n(\x) > h} \x^\square.
$$
\label{clusters}
\end{definition}
A cluster as a function of $n$ has a limiting shape if, appropriately scaled, it tends to a certain shape as $n \to \infty$, in some sense to be specified later.
By the model definition, all clusters are path connected.
Observe that from Definitions \ref{topplingdef} and \ref{clusters}, it follows that
\beq
\T_n \subset \V_n \subseteq \T_n \cup \partial\T_n,
\label{randvanTn}
\eeq
with $\partial\T_n$ the exterior boundary of $\T_n$.
We also define the ``lattice ball'':
\beq
\B_n = \bigcup_{i=1}^n \x_i^\square,
\label{ball}
\eeq
where the lattice sites $\x_1, \x_2, \ldots$ of $\mathbb{Z}^d$ are ordered in such a way that the Euclidean distance from the origin is non-decreasing.

\subsection{The rotor-router (RR,$h$) model}
\label{rotorrouterdef}

For the rotor router model, $c = 1$, so that in each toppling only one particle moves to a neighbor. Therefore, $h_{\max}=0$, and the model is defined for $h < 0$ only.
In words, every site holds the first $|h|$ particles that it receives, and sends every next particle to a neighbor, choosing the neighbors in a cyclic order. This means that after $2d$ topplings, every neighbor received a particle. Instead of a direction function with numerical values, we can think of an arrow being present at every site. In a toppling, the arrow is rotated to a new direction, and the particle is sent in this new direction.

Propp, Levine and Peres studied the case $h=-1$.
Levine and Peres have proven that for $h=-1$, the limiting shape of the rotor router model is a Euclidean sphere. More precisely, they showed (\cite{peres}, Thm. 1.1):
\beq
\lim_{n \to \infty}\lambda(n^{-1/d}\V_n \bigtriangleup B) = 0,
\label{rrresult}
\eeq
where $\lambda$ denotes $d$-dimensional Lebesgue measure, $\bigtriangleup$ denotes symmetric difference, and $B$ is the Euclidean sphere of unit volume centered at the origin in $\mathbb{R}^d$. This result is independent of $D_0$, and of any assignment of different unit vectors to possible values of $D(\x)$. Note that the scaling function is necessarily $n^{-1/d}$, since $|\V_n|=n$ by model definition.

In \cite{kleber}, a picture of $D_n$ is given for $h=-1$ and $n =$ 3 million and $D_0$ constant. The picture shows a circular shape with an intriguing seemingly self-similar pattern. Curiously enough, this picture actually shows the shape of $\T_n$ rather than $\V_n$, indicating that the limiting shape of $\T_n$ is also a sphere. This would be a stronger statement than (\ref{rrresult}), but it remains as yet unproven.

It has been noted that the shape of $\V_n$ for the rotor-router model is remarkably circular, that is, as close to a circle as a lattice set can get, for every $n$. However, the shape of $\T_n$ has not been studied before. We programmed the model for several values of $h$ and $n$, and observed for all these values that $\V_n \setminus \T_n$ is concentrated on the inner boundary of $\V_n$.

Our first main result of the rotor router model is the generalization of (\ref{rrresult}) to all $h \leq -1$, stated in Theorem \ref{RRshape}.
The proof in fact uses (\ref{rrresult}) as main ingredient.

By an entirely different method, in fact as a corollary of Theorem \ref{topplingbounds}, we moreover obtain the result that the rotor router shape contains a diamond of radius proportional to $(\frac n{2d-1-h})^{1/d}$. It is surprising that, in spite of the limiting shape of the RR,$h$ model being known, this is still a new result. But, in the words of Levine and Peres, the convergence in (\ref{rrresult}) does not preclude the formation of, e.g., holes close to the origin, as long as their volume is negligible compared to $n$. The only comparable previous result is that for $d=2$ and $h=-1$, the particle cluster contains a disk of radius proportional to $n^{1/4}$ \cite{levine}.

\subsection{The abelian sandpile (SP,$h$) model}
\label{sandpiledef}

For the abelian sandpile model, $c=2d$. Therefore, $h_{\max} = 2d-1$, and the model is defined for $h < 2d-1$. In each toppling, one particle moves to every neighbor, so that in fact the value of $D(\x)$ is irrelevant.

It follows that for the SP model, we can specify (\ref{randvanTn}) to
\beq
\V_n = \T_n \cup \partial \T_n.
\label{SPrandvanTn}
\eeq
The SP model as a growth model has received some attention in the cases $h=-1$ (``greedy'' sandpile) and $h=0$ (``non-greedy''), for which pictures of $H_n$ can be found \cite{kleber}. It is noted that the shape does not seem to be circular.
In Figure \ref{sandpileshapes}, we show a family of sandpile pictures for a range of values of $h$, obtained by programming the model in Matlab. We see the number of symmetry axes increasing as $h$ decreases. The shape appears to  become more circular as $h$ decreases. The shape for $h=2$ is observed to be a square, but for the other values of $h$ it seems to tend to a more complicated shape.

Again, we find that $\V_n \setminus \T_n$ is concentrated on the inner boundary of $\V_n$, for all observed values of $h$.

\begin{figure}[ht]
% \centerline{\includegraphics[width=13cm]{sandpileshapes}}
 \caption{$H_n$ for the sandpile model with different $h$ values. The number of particles ranges from $n = $ 60,000 ($h=2$) to 1,000,000 ($h=-20$).The gray-scale colors are such that a lighter color corresponds to a higher value of $H(\x)$.}
 \label{sandpileshapes}
\end{figure}

For the sandpile model, we have the following results.
Theorem \ref{cube} states that indeed the toppling cluster for the SP,$2d-2$ model is a ($d$-dimensional) cube, and the particle cluster tends to a cube as $n \to \infty$.

Theorem \ref{hminoneindig} states that for $h \to -\infty$, the limiting shape is a sphere.

Finally, we generalize some bounds for the scaling function that have been obtained by Le Borgne and Rossin \cite{rossin} in the case $d=2$, $0 \leq h \leq 2$, to all $d$ and $h$. This result is formulated in Theorem \ref{topplingbounds}. Our proof moreover allows to deduce that $\V_n$ for the sandpile model is simply connected for all $n$.

\subsection{The double router (DR,$h$) model}

In the double router model, $c=2$. Therefore, $h_{\max}=1$, and the model is defined for $h \leq 0$. In a toppling of this model, two particles are sent out in two different directions, such that after $d$ topplings, every neighbor received a particle.

Many variants of this model are possible, e.g., for $d=3$ one could choose $c=3$, which amounts to dividing the 6 neighbors in 2 groups of 3. However, to avoid confusion, we only discuss the above explained variant.

\begin{figure}[h]
% \centerline{\includegraphics[width=15cm]{mannashapes}}
 \caption{$H_n$ for the DR model with $h=0$ ($n =$ 110,000), $h=-1$ ($n =$ 100,000), and $h=-5$ ($n = $ 400,000), and $D_0(\x)=0$ for all $\x$. Sites with height 1 are colored white, sites with height 0 or negative are colored black.}
 \label{mannashapes}
\end{figure}

Figure \ref{mannashapes} shows $H_n$ for $h=0$, $d=2$ and $n = 110.000$. We ordered the unit vectors as $\e_0 =$ left, $\e_1 =$ right, $\e_2 =$ up, $\e_3 =$ down. Initially, $D(\x)=0$ for all $\x$. Figure \ref{mannaD} shows $D_n$ for the same case, to indicate that for this model we find complex patterns both for $H_n$ and $D_n$.

\begin{figure}[h]
% \centerline{\includegraphics[width=9cm]{mannaD}}
 \caption{$D_n$ for the DR model with $h=0$ ($n =$ 110,000), and $D_0(\x)=0$ for all $\x$. White corresponds to $D(\x)=2$, black to $D(\x)=0$.}
 \label{mannaD}
\end{figure}

As for the sandpile model, we see a symmetric, yet not circular shape, with notable flat edges. In Figure \ref{mannashapes}, also some other values of $h$ are shown, again indicating that the shape seems to become more and more circular with decreasing $h$. Indeed, Theorem \ref{hminoneindig} states this fact for both the sandpile model and the double router model.

Again, we observe that $\V_n \setminus \T_n$ is concentrated on the inner boundary of $\V_n$, for all observed values of $h$.

\section{Comparing the models}
\label{comparemodels}

From the model descriptions, it is clear that one SP toppling always equals $d$ DR topplings, as well as $2d$ RR topplings. Furthermore, our toppling definition (Definition \ref{topplingdef}) ensures that one DR toppling always equals two consecutive RR topplings. We exploit these relations to compare the clusters for different models.
We will with a subindex RR, DR or SP indicate the model that was used to obtain the growth cluster, and also use a subindex $h$ because we compare different $h$ values.

\begin{proposition}
For every fixed $h$ (for which the model is defined), $d$ and $n$,
\begin{enumerate}
\item
$\V_{n,SP,h} \subseteq \V_{n,DR,h} \subseteq \V_{n,RR,h}$,
\item
For all models $i = $ RR, DR and SP: $\V_{n,i,h-1} \subseteq \V_{n,i,h}$,
\item
$\V_{n,RR,h-1} \subseteq \V_{n,DR,h}$, and $\V_{n,RR,h-(2d-1)} \subseteq \V_{n,SP,h}$.
\end{enumerate}
\label{subsets}
\end{proposition}

\noindent
{\it Proof.}
The proof makes use of the abelian property of all models, that is, the property that the stabilized configuration $\eta_n$ does not depend on the order of topplings. We are therefore free to choose a convenient order.
Furthermore, a site $\x$ that at some instant during stabilization has either $T(\x)>0$, or $T(\x)=0$ and $H(\x)>h$, must belong to $\V_n$ in the final configuration, because by further topplings either $H(\x)$ or $T(\x)$, or both, increase.

{\it part 1.}
The initial configurations for all these three models are the same.
We first compare $\V_{n,SP,h}$ with $\V_{n,DR,h}$.
We choose, for both the SP and the DR model, to first perform all legal SP-topplings. Since every SP-toppling consists of $d$ DR-topplings, these are legal topplings for both models.
The configuration is now stabilized for the SP model, but there can be sites in $\V_{n,SP,h}$ that are not stable in the DR model, since $h_{max}$ is $2d-1$ for the SP model and 1 for the DR model. Therefore, in the MD model, possibly more topplings follow, so that $\V_{n,SP,h} \subseteq \V_{n,DR,h}$.

We compare $\V_{n,DR,h}$ with $\V_{n,RR,h}$ in an analogous manner, this time using that $h_{max}$ is 1 for the DR model and 0 for the RR,$h$ model, to get $\V_{n,DR,h} \subseteq \V_{n,RR,h}$.

{\it part 2.}
We start with comparing the RR,$h-1$ model with the RR,$h$ model.
We choose, for both initial configurations, to first perform all RR-topplings that stabilize the initial configuration of the RR,$h-1$ model. These are legal topplings for both configurations. The configuration is now stable for the RR,$h-1$ model, but all sites that have $H(\x)=0$ for this model, have $H(\x)=1$ for the RR,$h$ model. Therefore, for this model more topplings would follow, so that $\V_{n,RR,h-1} \subseteq \V_{n,RR,h}$.
The same reasoning can be applied for the DR and SP model.

{\it part 3.}
We first compare $\V_{n,DR,h}$ with $\V_{n,RR,h-1}$.
We choose, for both models, to first perform all DR-topplings that would stabilize the initial configuration with height $h-1$. These are legal topplings for both models. Then for both models, more topplings are needed. For the RR,$h-1$ model there can be sites with $H(\x)=1$, that are unstable. For the DR,$h$ model, the same sites have $H(\x)=2$, and therefore are also unstable. But since one DR-toppling equals two RR-topplings, the set of sites where the configuration changes by the extra topplings for the RR,$h-1$ model, is a subset of those for the DR,$h$ model.
Therefore, $\V_{n,RR,h-1} \subseteq \V_{n,DR,h}$.
The argument for the sandpile model is similar.
\qed

\section{Limiting shape results}
\label{shaperesults}

\subsection{The sandpile model in $\mathbb{Z}^d$, with $h=2d-2$}
\label{square}

As noted in Section \ref{definitionssection}, Figure \ref{sandpileshapes} indicates that the limiting shape for the sandpile model with $h=2$ and $d=2$ is a square. This section contains the proof of a more general statement for arbitrary dimension, that is, we prove that indeed the toppling cluster for the SP,$2d-2$ model is a ($d$-dimensional) cube, and the particle cluster tends to a cube as $n \to \infty$. However, we have no explicit expression for the scaling function $f(n)$, so that the scaled clusters $f(n)\V_n$ and $f(n)\T_n$ would tend to the unit cube, just that this scaling satisfies $n^{-1} \leq f(n) \leq \frac 12 n^{-1/d}- \frac 32$. Based on the calculations, we believe that $f(n)$ is $O(n^{-1/d})$.

\begin{theorem}
Let $\mathfrak{C}(r)$ be the cube $\bigcup_{\x: \max_i |x_i| \leq r}\x^\square$.
For every $n$ in the d-dimensional SP,$2d-2$ model, there is an $r_n$ such that
$$
\T_n = \mathfrak{C}(r_n),
$$
and
$$
\V_n = \mathfrak{C}(r_n) \cup \partial \mathfrak{C}(r_n).
$$
For all $n$, this $r_n$ satisfies
$$
\frac 12 n^{1/d} - \frac 32 \leq r_n \leq n.
$$
\label{cube}
\end{theorem}

We start by outlining the case $d=2$ as an example, for the sake of clarity. When simulating the model for several small values of $n$ - this can even be done by hand - one notices that the configuration is always such that it contains a central square with all boundary sites full, except for the corner sites, which have height $h_{max}-1$.  Outside this square, all sites have height $h_{max}-1$. We will call such a rectangular boundary a critical boundary. This is all the information we need to make an inductive argument in $n$.
Suppose, $\eta_n$ has a critical boundary, we add one grain to the origin and in the course of stabilization (to obtain $\eta_{n+1}$), one boundary site topples. As a consequence, all neighboring boundary sites topple, because they are all full and all in turn receive a grain. The two adjacent corner sites receive one grain, therefore become full. Therefore, after these topplings a new rectangular critical boundary is created. If any more boundary sites topple, we can reiterate this argument.
We conclude that the presence of a rectangular critical boundary is stable under additions inside this rectangle. In the special case where additions are made only to the origin, we conclude by symmetry that the shape of the critical boundary will always be square.

Below, we schematically show the argument for $n=7$.

\bigskip
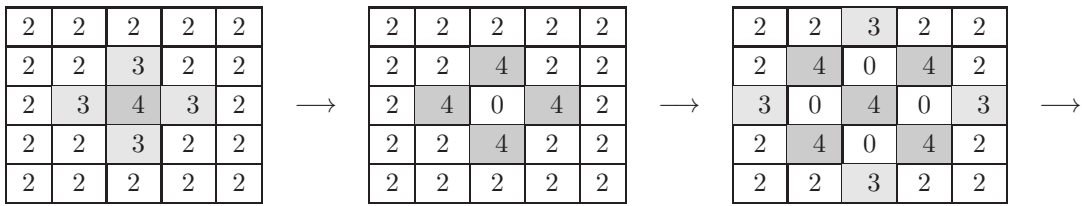
\begin{center}
\begin{tabular}{cccccc}
\begin{tabular}{|c|c|c|c|c|}
\hline
 2 & 2 & 2 & 2 & 2\\
\hline
 2 & 2 & \cellcolor{lichtgrijs} 3 & 2 & 2\\
\hline
 2 & \cellcolor{lichtgrijs} 3 & \cellcolor{grijs} 4 & \cellcolor{lichtgrijs} 3 & 2\\
\hline
 2 & 2 & \cellcolor{lichtgrijs} 3 & 2 & 2\\
\hline
 2 & 2 & 2 & 2 & 2\\
\hline
\end{tabular} & $\longrightarrow$ &
\begin{tabular}{|c|c|c|c|c|}
\hline
 2 & 2 & 2 & 2 & 2\\
\hline
 2 & 2 & \cellcolor{grijs} 4 & 2 & 2\\
\hline
 2 & \cellcolor{grijs} 4 & 0 & \cellcolor{grijs} 4 & 2\\
\hline
 2 & 2 & \cellcolor{grijs} 4 & 2 & 2\\
\hline
 2 & 2 & 2 & 2 & 2\\
\hline
\end{tabular} & $\longrightarrow$ &
\begin{tabular}{|c|c|c|c|c|}
\hline
 2 & 2 & \cellcolor{lichtgrijs} 3 & 2 & 2\\
\hline
 2 & \cellcolor{grijs} 4 & 0 & \cellcolor{grijs} 4 & 2\\
\hline
 \cellcolor{lichtgrijs} 3 & 0 & \cellcolor{grijs} 4 & 0 & \cellcolor{lichtgrijs} 3\\
\hline
 2 & \cellcolor{grijs} 4 & 0 & \cellcolor{grijs} 4 & 2\\
\hline
 2 & 2 & \cellcolor{lichtgrijs} 3 & 2 & 2\\
\hline
\end{tabular} & $\longrightarrow$ \\
\end{tabular}
\end{center}

\medskip
\noindent
{\it We start with $\eta_7$ plus one extra grain at the origin. Full boundary sites are colored lightgrey, unstable sites grey. First the origin topples, causing the full boundary sites to become unstable. Next, these boundary sites topple, causing the corner sites to be unstable. The first full sites of the new boundary are created.}
\medskip

\begin{center}
\begin{tabular}{ccccc}
\begin{tabular}{|c|c|c|c|c|}
\hline
 2 & \cellcolor{lichtgrijs} 3 & \cellcolor{lichtgrijs} 3 & \cellcolor{lichtgrijs} 3 & 2\\
\hline
 \cellcolor{lichtgrijs} 3 & 0 & 2 & 0 & \cellcolor{lichtgrijs} 3\\
\hline
 \cellcolor{lichtgrijs} 3 & 2 & \cellcolor{grijs} 4 & 2 & \cellcolor{lichtgrijs} 3\\
\hline
 \cellcolor{lichtgrijs} 3 & 0 & 2 & 0 & \cellcolor{lichtgrijs} 3\\
\hline
 2 & \cellcolor{lichtgrijs} 3 & \cellcolor{lichtgrijs} 3 & \cellcolor{lichtgrijs} 3 & 2\\
\hline
\end{tabular} & $\longrightarrow$ &
\begin{tabular}{|c|c|c|c|c|}
\hline
 2 & \cellcolor{lichtgrijs} 3 & \cellcolor{lichtgrijs} 3 & \cellcolor{lichtgrijs} 3 & 2\\
\hline
 \cellcolor{lichtgrijs} 3 & 0 & 3 & 0 & \cellcolor{lichtgrijs} 3\\
\hline
 \cellcolor{lichtgrijs} 3 & 3 & 0 & 3 & \cellcolor{lichtgrijs} 3\\
\hline
 \cellcolor{lichtgrijs} 3 & 0 & 3 & 0 & \cellcolor{lichtgrijs} 3\\
\hline
 2 & \cellcolor{lichtgrijs} 3 & \cellcolor{lichtgrijs} 3 & \cellcolor{lichtgrijs} 3 & 2\\
\hline
\end{tabular}\\
\end{tabular}
\end{center}
\medskip
\noindent
{\it Next, the corner sites topple. The configuration now has a new square critical boundary. A further toppling of the origin, after which we obtain $\eta_8$, does not change this new boundary anymore.}
\medskip

In the proof below, we generalize this critical boundary to arbitrary $d$ and make the details more precise, but the idea will remain the same.

In the course of the proof, we will need the following lemma.

\begin{lemma}
Let $\phi_r$ be the following configuration: $H(\0) = 2d$, and for all other $\x \in \mathfrak{C}_r$, $H(\x) = h_{max}$, otherwise $H(\x) = h_{max}-1$.
If the configuration $\phi_r$ is stabilized, then during stabilization, every site $\x$ topples exactly $\delta_{\x,r} = \max \{r+1-\max_i|x_i| ,0\}$ times.
\label{squaretopplelemma}
\end{lemma}

\noindent{\it Proof.}
We choose to order the topplings into \emph{waves} \cite{priezzhev}, that is, in each wave, we topple the origin once and then all other sites that become unstable, except the origin again. No site can topple more than once during a wave.

In the first wave, all sites in $\mathfrak{C}_r$ topple once, because it is maximally filled. The wave stops at the boundary, because the sites outside $\mathfrak{C}_r$ have at most one toppling neighbor, therefore they can not become unstable. All sites in $\mathfrak{C}_{r-1}$ then have $2d$ once-toppling neighbors, so their particle number does not change. The sites in $\mathfrak{C}_r \setminus \mathfrak{C}_{r-1}$ have less toppling neighbors, so their particle number becomes at most $2d-2$.

Therefore, the effect of the first wave on $\phi_r$, is to make it, for all sites in $\mathfrak{C}_r$, at most equal to $\phi_{r-1}$, with equality for all sites with $\mathfrak{C}_{r-1}$. It follows that the next wave will topple all sites in  $\mathfrak{C}_{r-1}$ once. Continuing this argument, the result stated in the lemma follows.
\qed

\noindent
{\it Proof of Theorem \ref{cube}.}
We will use induction in $n$. For that, we choose to obtain the final configuration as follows: the $n$ particles are added to the origin one by one, each time first stabilizing the current configuration through topplings. Due to abelianness, this procedure will give the same final configuration as when all $n$ particles are added simultaneously, before toppling starts.
We will show that during this procedure, only configurations in $\C_r$ are encountered, where $\C_r$ is the set of configurations that are as follows:
\begin{enumerate}
\item
$\T_n = \mathfrak{C}_r$,
\item
For all $\x \in \T_n$ such that $\max_i|x_i|=r$, $T_n(\x)=1$,
\item
$\V_n = \mathfrak{C}_r \cup \partial \mathfrak{C}_r$,
\item
For all $\x \in \V_n\setminus \T_n$, $H_n(\x)=h_{max}=2d-1$,
\item
For all $\x \not \in \V_n$, $H_n(\x) = 2d-2$.
\end{enumerate}
In words, a configuration in $\C_r$ is such that $\T_n$ is a cube, where all inner boundary sites toppled exactly once, and all outer boundary sites, i.e., sites that do not belong to $\T_n$ but have a neighbor in $\T_n$, have particle number $h_{max}$. Thus, $\T_n$ is surrounded by $2d$ maximally filled square ``slabs'' of size $(2r+1)^{d-1}$. All sites outside $\V_n$ have particle number $2d-2$, by model definition.

After the first particle is added, $H_1(\0) = h_{max}$, and no site has toppled yet. After the second particle is added, the origin topples once, and all neighbors of the origin receive one particle. Therefore, $\eta_2 \in \C_0$.

We will now show that if we assume $\eta_n \in C_r$ for some $r$, then either $\eta_{n+1} \in C_r$ or $\eta_{n+1} \in C_{r+1}$. To prove this, we will add a particle at the origin to configuration $\eta_n$ and start stabilizing. First we show that if during stabilization the configuration leaves $\C_r$, then it enters $\C_{r+1}$. Then we show that in this last case, it does not leave $\C_{r+1}$ during further stabilization.

If no site $\x \in \V_n$ such that $\max_i |x_i|=r$ topples during stabilization, then $\eta_{n+1}$ remains in $C_r$.
But if one such site topples, then one site in $\V_n \setminus \T_n$ becomes unstable, and also topples. This site is in one of the $2d$ maximally filled slabs. If one site of such a slab topples, then the entire slab must topple, because all of its sites will in turn receive a particle from a toppling neighbor.
After the entire slab toppled once, all neighbors of the slab received a particle, and the configuration can be described as follows: The toppled cluster is a rectangle of size $(2r+1)^{d-1}(2r+2)$, surrounded by maximally filled slabs, two of them cubic and the rest rectangular.

But due to symmetry, if one slab topples once then this must happen for all slabs. We can choose in what order to topple the slabs; suppose we divide them in $d$ opposite pairs. After we toppled the first pair of slabs, the configuration is a central rectangle of size $(2r+1)^{d-1}(2r+3)$, surrounded by maximally filled slabs, two of them cubic and the rest rectangular. After we toppled the $k^{\rm th}$ pair of (meanwhile possibly rectangular) slabs, the configuration is a central rectangle of size $(2r+1)^{d-k}(2r+3)^k$, surrounded by maximally filled slabs, so that after all slabs toppled, the toppled cluster is again a cube, now of size $(2r+3)^d$, centered at the origin and surrounded by $2d$ maximally filled square slabs.
In other words, after these topplings the configuration is in $\C_{r+1}$. It follows that if $\eta_n \in \C_r$ for some $r$, then $\eta_{n+1} \in \C_{r'}$, for some $r'\geq r$. It now remains to show that $r'$ can only have the values $r$ or $r+1$.

We use lemma \ref{squaretopplelemma}.
A configuration $\eta \in \C_r$, plus an addition at the origin, has at every site at most the number of particles given by $\phi_{r+1}$. Therefore, if we suppose $\eta_n \in \C_r$, then upon addition of a particle at the origin, by abelianness, the number of topplings for every $\x$ will at most be $\delta_{\x,r+1}$. In particular, the outer boundary sites of $\T_n$ will topple at most once.

The induction proof is now completed, therefore we now know that, for all $n$, $\eta_n \in \C_{r}$ for some $r$. We have also shown that $r_n \leq n$. From the description of $\C_r$, we see that $\T_n$ is always a cube, and that $\V_n$ is more and more like a cube as $r_n$ increases.

To conclude that $\V_n$ tends to a cube as $n \to \infty$, we finally need to show that $r_n$ increases with $n$. For this, note that in $\eta_n$ every site needs to be stable. Since $h = h_{max}-1$, every site can accommodate at most one extra particle. Therefore, if $\eta_n \in C_r$, then $r \geq \frac 12 n^{1/d} - \frac 32$.
\qed

\subsection{The rotor-router model with $h<-1$}
\label{RRshapesection}

In this section, we will prove the following result:

\begin{theorem}
The limiting shape of the particle cluster for the rotor-router model is a sphere for every $h \leq -1$. More precisely,
$$
\lim_{n \to \infty} \lambda\left((\frac n{|h|})^{-1/d}\V_n \bigtriangleup B\right) = 0.
$$
\label{RRshape}
\end{theorem}

The proof will also reveal the following about the rotor-router toppling cluster:

\begin{corollary}
The toppling cluster of the RR,$h$ model contains a cluster $\W_n$, with
$$
\lim_{n \to \infty} \lambda\left((\frac n{|h|+1})^{-1/d}\W_n \bigtriangleup B\right) = 0.
$$
\label{RRtopplingcluster}
\end{corollary}

Theorem \ref{RRshape} has been proven for the case $h=-1$, by Levine and Peres, see (\ref{rrresult}). We will use their result to prove the theorem for other values of $h$.
The strategy of the proof will be as follows. First, we introduce a slightly different version of the RR,$h$ model, which we will call the $k$-color rotor router, or RR',$k$ model. We show that the limiting shape of this model is a sphere for all $k$, with scaling function $(\frac nk)^{-1/d}$. Informally, the $k$-color model can be viewed as $k$ iterations of the RR,-1 model. The $n$ particles at the origin are equally divided in $k$ different-colored groups, and for each color the particles spread out until there is at most one particle of each color at each site. From (\ref{rrresult}), we then know that each color region is approximately spherical, so we get a final configuration that looks like $k$ almost completely overlapping, approximately spherical, different-colored regions.

Then we show that $\V_{n,RR,h}$ and $\V_{n,RR',k}$, with $h = -k$, differ only in a number of sites that is $o(n)$.
This number is at most the number of sites in $\V_{n,RR',k}$ where not every color is present. In the proof of this second point, we will stabilize $\eta_0$ for the RR,$h$ model by first performing all topplings needed to stabilize the RR',$-h$ model. Some of these topplings may be illegal for the RR,$h$ model, so we will reverse them by performing untopplings. The main difficulty in the proof is to show that we do reach $\eta_n$ by this procedure. First we define untopplings.

\begin{definition}
An untoppling of site $\x$ in configuration $\eta$ consists of the following operations:
\begin{itemize}
\item
$T(\x) \rightarrow T(\x) - 1$,
\item
$H(\x) \rightarrow H(\x) + c$, with $c$ according to the model,
\item
$D(\x) \to (D(\x)-c) \mod 2d$,
\item
$H(\x + \e_i) \rightarrow H(\x + \e_i) - 1$, with $i = D(\x), (D(\x)+1) \mod 2d, \ldots, (D(\x)+c-1) \mod 2d$.
\end{itemize}
\end{definition}

We call an untoppling of site $\x$ {\em legal} if $T(\x)>0$ before the untoppling. We see that a legal untoppling of site $\x$ is precisely the undoing of a toppling of site $\x$, in the sense that the particles that were sent out of site $\x$ in the last toppling, return to this site, and the value of $D(\x)$ returns to its previous value.

The following proposition provides a more elaborate description of $\eta_n$ in the case where we allow illegal topplings, that we will need in the proof of Theorem \ref{RRshape} (we remark that the abelian property still holds for combinations of legal and illegal topplings). We introduce the notion of optimality to distinguish $\eta_n$ from configurations that are stable and allowed, but that cannot be obtained from $\eta_0$ by legal topplings. An example for the rotor router model is as follows. Suppose that, after reaching $\eta_n$ by legal topplings, there is somewhere a closed loop of sites with height 0 such that ``the arrows form a cycle'', i.e., if we topple all of them, they would all send and receive one particle.
After these topplings (of which at least the first would be illegal), the height function is still $H_n$, but the toppling number of these sites has increased by 1, so we obtained a different, stable and allowed configuration.
We remark that a set of sites forming such a ``cycle of arrows'' cannot be found in $\T_n$. Since all sites that toppled are full, in the last toppling in this set a particle must have left the set, so that there is an arrow pointing out of the set (see also \cite{walkers}, section III).
Therefore, if there is a set of sites forming a cycle of arrows, then at least one of these sites did not topple.

\begin{proposition}
Call a stable configuration {\em optimal} if every sequence of legal untopplings leads to an unstable configuration.
For every model, $\eta_n$ is the unique optimal, stable and allowed configuration that can be reached from $\eta_0$ by topplings, either legal or illegal, and legal untopplings.
\label{optimal}
\end{proposition}

\noindent{\it Proof.}
We defined $\eta_n$ before as the unique stable configuration that can be reached from $\eta_0$ by legal topplings. It follows from this definition that $\eta_n$ is allowed, since if sites can only topple when they are unstable, then we automatically obtain for all $\x$ where $T_n(\x)>0$, that $H_n(\x)\geq 0$.

To prove that $\eta_n$ is optimal, we proceed by contradiction. Suppose that, starting from $\eta_n$, there is a sequence of legal untopplings such that a stable configuration $\xi$ is obtained. Then $\xi$ can be obtained from $\eta_0$ by a sequence of legal or illegal topplings. Call $T'$ the toppling function for this sequence, then $0 \leq T'(\x) \leq T_n(\x)$ for all $\x$ (because we undid some topplings of $T_n$), but $T' \neq T_n$. By abelianness, $\xi$ depends only on $\eta_0$ and $T'$. Thus, we can choose the order of the topplings according to $T'$, to obtain $\xi$. There cannot be an order such that all topplings are legal, otherwise $\xi$ cannot be different from $\eta_n$. We will therefore choose the order such that first all possible legal topplings are performed, and then the rest (as an example, suppose $T'(\0)=0$. Then we must start with an illegal toppling, since in $\eta_0$ all sites but the origin are stable).
After all possible legal topplings according to $T'$, we have a configuration with at least one unstable site, since we did not yet reach $\eta_n$. In the remainder, this site is not toppled, because otherwise another legal toppling could have been added to the first legal topplings. Therefore, $\xi$ cannot be stable.

We so far established that $\eta_n$ is optimal, stable and allowed. Now we prove that it is unique by deriving another contradiction. Suppose there is an other optimal, stable and allowed configuration $\zeta$ that can be reached from $\eta_0$ by topplings and legal untopplings. Call $T''$ the toppling function for this sequence.
Then $T''(\x) \geq T_n(\x)$ for all $\x$. We can see this as follows: we choose the toppling order such that we first perform all topplings that are also in $T_n$. If there were some sites $\y$ where $T''(\y) < T_n(\y)$, then at this point at least one of these sites would be unstable. But in the remainder of $T''$, this site would not topple again, so in that case $\zeta$ could not be stable.

Call $\tau(\x) = T''(\x)-T_n(\x)$. Since $\zeta \neq \eta_n$, $T''\neq T_n$. Suppose we first perform the topplings according to $T_n$, and then according to $\tau$.
Then we have that $\zeta$ can be reached from $\eta_n$ by performing $\tau(\x)$ topplings for every site $\x$, and vice versa, that $\eta_n$ can be reached from $\zeta$ by performing $\tau(\x)$ legal untopplings for every site $\x$. But then $\zeta$ cannot be optimal.
\qed

\medskip
\noindent
{\em Proof of Theorem \ref{RRshape}}
The RR',$k$ model is defined as follows: As in the RR,$h$ model, $c=1$. Initially, there are $n$ particles at the origin; at every other site there are $-k$ particles. We choose $n$ a multiple of $k$, and divide the $n$ particles into $k$ groups of $n/k$ particles, each group with a different color. The height $H(\x)$ of a site $\x$ is now defined as $-k$ plus the total number of particles present at $\x$, of either color. However, now we call a site stable if it contains at most one particle of each color. If a particle arrives at a site where its color is already present, then that particle will be sent to a neighbor in the subsequent toppling. Note that in this model, a site $\x$ can be unstable even if $H(\x)<0$.
We furthermore restrict the order of legal topplings for this model.
We say that in this model, first all legal topplings for the first color should be performed, then for the next, etc.
Then for each color, the model behaves just as the $RR,-1$ model, with for each color a new initial $T$ and $D$. Therefore, with this toppling order the final configuration consisting of $H_n$, $T_n$ and $D_n$, is well-defined.

We now show that the limiting shape of the RR',$k$ model is a sphere.
We perform first all the topplings with particles of the first color, say, red. From \cite{peres}, Thm. 2.1 it follows that these particles will form a cluster $\V^1_{n/k}$ close to the lattice ball $\B_{n/k}$. In fact, the number of sites in $\V^1_{n/k} \bigtriangleup \B_{n/k}$ is $o(n)$. We will denote this as
$$
|\V^1_{n/k} \bigtriangleup \B_{n/k}| \leq f(n/k) = o(n).
$$
$D$ is now different from $D_0$. However, if next we stabilize for the blue particles, these will again form a cluster close to the lattice ball $\B_{n/k}$, since the result (\ref{rrresult}) does not depend on the initial $D$. Of course, this cluster $\V^2_{n/k}$ need not be the same as $\V^1_{n/k}$.

When all sites are stable, we have $\V_{n,RR',k} = \bigcup_{i=1}^k \V^i_{n/k}$. We also define a cluster $\W_{n,k} = \bigcap_{i=1}^k \V^i_{n/k}$. From the above, the number of sites in the difference of both these clusters with $\B_{n/k}$, is at most $k f(n/k)$. Thus
\beq
|\V_{n,RR',k} \bigtriangleup \B_{n/k}| \leq k f(n/k) = o(n).
\label{colorshape}
\eeq
We will call $\X_n = \V_{n,RR',k} \setminus \W_{n,k}$, so that $\X_n$ contains $|\X_n| \leq 2k f(n/k)$ sites. Sites in $\W_{n,k}$ contain $k$ particles, sites in $\X_n$ contain less than $k$ particles.

Now we compare this model with the RR,$h$ model, with $h=-k$. Disregarding the colors, the initial configuration is the same for both models. Suppose we perform in the RR,$h$ model all the topplings that are performed as above in the RR',$-h$ model. The configuration is then as follows: for all $\x \in \W_{n,k}$, $H(\x)=0$ and $T(\x) \geq 0$, and for all $\x' \in \X_n$, $-k \leq H(\x') < 0$ and $T(\x') \geq 0$.

This configuration is possibly not $\eta_n$ for the RR,$h$ model, since it is possibly not allowed, if there are sites $\x'$ in $\X_n$ with $T(\x')>0$. We can reach an allowed configuration by performing legal untopplings; it will appear that we then in fact reach $\eta_n$.

Suppose first that only one untoppling is needed. It might be that the neighbor $\y$ that returned a particle, now has less than $k$ particles, so that $H(\y)<0$, while $T(\y)>0$. In that case, $\y$ would now also have to untopple. We call this process an untoppling avalanche. The avalanche stops if an allowed (stable) configuration is encountered, that is, if the last particle came from a site that did not topple. An untoppling avalanche consists of untoppling neighbors, each one passing a particle to the previous one. Therefore, an untoppling avalanche of arbitrary length, changes the particle number of only two sites in the configuration.

In our case, at most $(k-1)|\X_n|$ untoppling avalanches are required. Each avalanche stops if the last particle came from a site that did not topple. This will change the configuration in at most $2(k-1)|\X_n|$ sites. After all these untopplings, the configuration is allowed and stable, i.e., we have for all sites with $T(\x)>0$, that $H(\x)=0$.
We can only perform legal untopplings at sites with $T(\x)>0$, so that we cannot perform a sequence of legal untopplings in a closed loop of sites. Therefore, any further sequence of legal untopplings would increase $H(\x)$ for at least one of these sites, and thus lead to an unstable configuration.

We have thus reached an optimal configuration, so we have reached $\eta_n$. Therefore, at this point we can conclude
$$
|\V_{n,RR,h} \bigtriangleup \V_{n,RR',k}| \leq  2(k-1)|\X_n| \leq 4k(k-1)f(n/k) = o(n).
$$
Combined with (\ref{colorshape}), this leads to Theorem \ref{RRshape}.
\qed

\medskip
\noindent
To prove the corollary, we define $\W_{n,h} = \bigcup_{\x:\x \in \V_{n,RR,h},H_n(x)=|h|} \x^\square$. By the above argument, we have, as for $\W_{n,k}$, that $|\W_{n,h} \bigtriangleup \B_{n/|h|}| = o(n)$.

Suppose, in the RR,$h$ model, we first perform all topplings that are needed to stabilize for the RR,$h-1$ model, as in the proof of Proposition \ref{subsets} part 2. The particle cluster then contains a cluster $\W_{n,h-1}$, with $\lim_{n \to \infty} \lambda\left((\frac n{|h-1|})^{-1/d}\W_n \bigtriangleup B\right) = 0$.

Because the sites in $\W_{n,h-1}$ contain $|h-1| = |h|+1$ particles, every site in $\W_{n,h-1}$ is unstable in the RR,$h$ model. Therefore, $\W_{n,h-1} \subseteq \T_{n,RR,h}$.
\qed

\subsection{The double router model, and the sandpile model with $h \to -\infty$}

In discussing Figures \ref{sandpileshapes} and \ref{mannashapes}, we observed that the DR and SP shapes seem to become more circular as $h$ decreases.
Proposition \ref{subsets} and Theorem \ref{RRshape} can indeed be combined to give the following result:
\begin{theorem}
The limiting shape of the SP,$h$ and the DR,$h$ model, for $h \to -\infty$, is a sphere. More precisely,
$$
\lim_{h \to -\infty} \limsup_{n \to \infty} \lambda\left((\frac n{|h|})^{-1/d}\V_n \bigtriangleup B\right) = 0.
$$
\label{hminoneindig}
\end{theorem}

\noindent{\it Proof.}
From Proposition \ref{subsets}, for all $h \leq -1$,
\beq
\V_{n,RR,h-(2d-1)} \subseteq \V_{n,SP,h} \subseteq \V_{n,RR,h},
\label{SPincl}
\eeq
and
\beq
\V_{n,RR,h-1} \subseteq \V_{n,DR,h} \subseteq \V_{n,RR,h}.
\label{DMincl}
\eeq
We will discuss the DR,$h$ model first.

From Theorem \ref{RRshape}, we know that the particle cluster of the RR,$h$ model, scaled by $(\frac n{|h|} )^{-1/d}$, tends for every fixed value of $h$ to the unit volume sphere as $n \to \infty$, i.e.,
$$
\lim_{n \to \infty} \lambda \left( (\frac n{|h|} )^{-1/d} \V_{n,RR,h} \bigtriangleup B \right) = 0,
$$
and thus we also know (note that since $h$ is negative, $|h-1| = |h|+1$)
$$
\lim_{n \to \infty} \lambda \left( (\frac n{|h|} )^{-1/d} \V_{n,RR,h-1} \bigtriangleup B \right) =
\lambda\left( \left[(\frac {|h|}{|h|+1} )^{1/d}B \right] \bigtriangleup B\right) = 1-\frac {|h|}{|h|+1}.
$$

Because of (\ref{DMincl}), it follows that
$$
0 \leq \limsup_{n \to \infty} \lambda \left( (\frac n{|h|} )^{-1/d} \V_{n,DR,h} \bigtriangleup B \right) \leq 1-\frac {|h|}{|h|+1},
$$
so that
$$
\lim_{h \to -\infty}\limsup_{n \to \infty} \lambda \left( (\frac n{|h|})^{-1/d} \V_{n,DR,h} \triangle B \right) = 0.
$$
The argument is similar for the SP,$h$ model.
\qed

\subsection{The toppling cluster for the sandpile model}
\label{topplingcluster}

In this subsection, we generalize Theorem 4 of Le Borgne and Rossin \cite{rossin}, who studied the particle cluster of the SP,$h$ model for $d=2$ and $h = 0,1$ and $2$, to arbitrary $d$ and $h$.

\begin{theorem}
Let $\mathfrak{D}(r)$ be the diamond $\bigcup_{\x: \sum_{i=1}^d |x_i| \leq r}\x^\square$, and $\mathfrak{C}(r)$ the cube $\bigcup_{\x: \max_i |x_i| \leq r}\x^\square$.
\begin{enumerate}
\item
In the SP,$h$ model, for every $n$ and $h \leq 2d-2$, there is an $r_n$ such that
$$
\mathfrak{D}(r_n-1) \subseteq \T_n \subseteq \mathfrak{C}(r_n).
$$
and
$$
\mathfrak{D}(r_n) \subseteq \V_n \subseteq \mathfrak{C}(r_n+1).
$$
\item
This $r_n$ satisfies
$$
\begin{array}{ll}
\left(\frac n{2d-1-h}\right)^{1/d}-3  \leq 2r_n & \mbox{for all $n$ and~} h \leq 2d-2, \\
2r_n \leq \left(\frac {dn}{d-h}\right)^{1/d} + o(n^{1/d}) & \mbox{for $n \to \infty$ and~} h < d.
\end{array}
$$
\end{enumerate}
\label{topplingbounds}
\end{theorem}

\begin{corollary}
$\V_{n,RR,h}$ contains the diamond $\mathfrak{D}\left(\frac 12 \left(\frac n{2d-1-h}\right)^{1/d}-\frac 32 \right)$ for all $h<0$, and is contained in the cube $\mathfrak{C}(\frac 12 (\frac {dn}{d-h})^{1/d})$ for all $h<-d$.
\label{rrcorollary}
\end{corollary}
The corollary follows from combining the theorem with Proposition \ref{subsets} part 1 and 3.
\begin{remark}
Note that the first inequality for $r_n$ agrees with that in Theorem \ref{cube}, for the sandpile model with $h=2d-2$. We can not use the second inequality in this case, because that is only valid for $h<d$. From combining Theorem \ref{cube} with Proposition \ref{subsets} part 3, we find that for $-d\leq h<0$, $\V_{n,RR,h}$ is contained in the cube $\mathfrak{C}(n+1)$.
\end{remark}

\noindent{\it Proof of Theorem \ref{topplingbounds}, part (1)}
We prove the inequality for $\T_n$. The inequality for $\V_n$ then follows from (\ref{SPrandvanTn}).

We follow the method of Le Borgne and Rossin, who introduce the following stabilization procedure: Starting at $t=0$ in the initial configuration with $n$ particles at the origin, each time step every unstable site of the current configuration is stabilized. For example, to obtain $\eta^1$ we topple only the origin, as many times as legally possible. At $t=2$, we topple only the neighbors of the origin, et cetera. We will call this the Le Borgne-Rossin procedure. Thus, $\eta^t$ contains stable sites that toppled at time $t$, and unstable sites that received grains. $T^t$ contains for every site the total number of topplings up to time $t$.
For every site $\x$, we then have
\beq
T^{t+1}(\x) = \max \{\left\lfloor \frac 1{2d}\left(H^0(\x) + \sum_{\y\sim \x}T^t(\y)\right) \right\rfloor,0\},
\label{topplings}
\eeq
where $\y\sim \x$ means that $\y$ is a neighbor of $\x$. We remark that this formula should correspond to formula (2) of \cite{rossin}, but, apparently due to a misprint, there the term $H^0(\x)$ was omitted.

If we divide the grid $\mathbb{Z}^d$ in two (``checkered'') subgrids $\{\x: \sum_{i=1}^d x_i ~\mbox{is even} \}$ and $\{\x: \sum_{i=1}^d x_i ~\mbox{is odd} \}$, then at every time $t$, at most one subgrid contains unstable sites, because the neighbors of sites of one subgrid are all in the other subgrid, and we started with only one unstable site.
We will consider next-nearest neighbor pairs $\x$ and $\z$, i.e., $\x$ and $\z$ are in the same subgrid, at most 2 coordinates differ, and $\sum_i \left||x_i|-|z_i|\right| = 2$, or in other words, $\z$ is a neighbor of a neighbor of $\x$, $\x$ itself being excluded.

First, we will prove by induction in $t$ that for every next-nearest neighbor pair $\x$ and $\z$, with $d(\x) \leq d(\z)$, where $d(\x)$ is the Euclidean distance of $\x$ to the origin, $T^t(\x) \geq T^t(\z)$. We remark that for a next-nearest neighbor pair $\x'$ and $\z'$ with $d(\x') = d(\z')$, we must have that $\x'$ is equal to $\z'$ up to a permutation of coordinates. The model should be invariant under such a permutation, therefore, in that case we have $T^t(\x') = T^t(\z')$.

The statement is true at $t=0$, since at $t=0$, we have $T^0(\x)=0$ for all $\x$. It is also true at $t=1$, since at $t=1$, only the origin topples.

Now suppose that the statement is true at time $t$, and let $\x_1$ and $\x_2$ be two next-nearest neighbor sites that do not topple at time $t$.
Suppose $d(\x_1) < d(\x_2)$. Then by assumption, $\sum_{\y_1\sim \x_1}T^t(\y_1) \geq \sum_{\y_2\sim \x_2}T^t(\y_2)$, because the neighbors of $\x_1$ and $\x_2$ can be grouped into next-nearest neighbor pairs with $d(\y_1) < d(\y_2)$.
Furthermore, we always have $H^0(\x_1) \geq H^0(\x_2)$, since $H^0(\0) = n$, and $H^0(\x) = h$ if $\x \neq \0$.

Inserting this in (\ref{topplings}), we obtain $T^{t+1}(\x_1) \geq T^{t+1}(\x_2)$, so that the statement remains true at $t+1$ for one subgrid.
The other subgrid must contain only sites that do not topple at $t+1$, so that for all $\x$ in this subgrid, $T^{t+1}(\x) = T^t(\x)$. Therefore, the statement remains true at $t+1$ for both subgrids.

Now let $\x$ be a site with maximal $r = \max_i |x_i|$ where $T(\x)>0$ (by symmetry, $\x$ cannot be unique). Then no sites outside $\mathfrak{C}(r)$ can have toppled. Furthermore, all next-nearest neighbors of $\x$ that are closer to the origin have also toppled, and subsequently all next-nearest neighbors of those next-nearest neighbors, that are still closer to the origin, etc.
Then we use symmetry to find that all sites in $\mathfrak{D}(r)$ and in the same subgrid as $\x$ must also have toppled. An example for $d=2$ is depicted in Figure \ref{diamondexample}.

\bigskip

\begin{center}
\begin{tabular}{c|c|c|c|c|c|c}
$\cdot$ & $\cdot$ & $z$ & \cellcolor{lichtgrijs} $\cdot$ & $z$ & $\cdot$ & $\cdot$ \\
\hline
$\cdot$ & $y$ & \cellcolor{lichtgrijs} $\cdot$ & \cellcolor{lichtgrijs} $z$ & \cellcolor{lichtgrijs} $\cdot$ & $z$ & $\cdot$ \\
\hline
$\bf{x}$ & \cellcolor{lichtgrijs} $\cdot$ & \cellcolor{lichtgrijs} $y$ & \cellcolor{lichtgrijs} $\cdot$ & \cellcolor{lichtgrijs} $z$ & \cellcolor{lichtgrijs} $\cdot$ & $z$ \\
\hline
\cellcolor{lichtgrijs} $\cdot$ & \cellcolor{lichtgrijs} $y$ & \cellcolor{lichtgrijs} $\cdot$ & \cellcolor{lichtgrijs} $\mathbf{0}$ & \cellcolor{lichtgrijs} $\cdot$ & \cellcolor{lichtgrijs} $z$ & \cellcolor{lichtgrijs} $\cdot$\\
\hline
$z$ & \cellcolor{lichtgrijs} $\cdot$ & \cellcolor{lichtgrijs} $z$ & \cellcolor{lichtgrijs} $\cdot$ & \cellcolor{lichtgrijs} $z$ & \cellcolor{lichtgrijs} $\cdot$ & $z$\\
\hline
$\cdot$ & $z$ & \cellcolor{lichtgrijs} $\cdot$ & \cellcolor{lichtgrijs} $z$ & \cellcolor{lichtgrijs} $\cdot$ & $z$ & $\cdot$ \\
\hline
$\cdot$ & $\cdot$ & $z$ & \cellcolor{lichtgrijs} $\cdot$ & $z$ & $\cdot$ & $\cdot$\\
\end{tabular}
\end{center}
\medskip

{\it Figure \ref{diamondexample}. Let $x$ be a site with maximal $r = \max_i|x_i|$ that toppled. In this example, $r=3$. By the fact that all next-nearest neighbors that are close to the origin have also toppled, we find that sites $y$, and the origin, have also toppled. By symmetry, we find that sites $z$ have also toppled. Then we conclude that all sites in $\mathfrak{D}(3)$ (colored lightgrey) and in the same subgrid as $x$ have toppled, and no sites outside $\mathfrak{C}(3)$ in the same subgrid can have toppled.}
\label{diamondexample}

\bigskip

Furthermore, for at least one of the neighbors $\y$ of $\x$, we must have $T(\y)>0$, because $\T_n$ is path connected. This neighbor has $\max_i|y_i| \geq r-1$. For this neighbor, we can make the same observation. It follows that all sites in $\mathfrak{D}(r-1)$ have toppled, so that the first part of the theorem follows.
\qed

\medskip\noindent
To prove the inequalities for $r_n$, we need the following lemma:
\begin{lemma}
Let $\rho_n$ be the average value of $H_n(\x)$ in $\T_n$, $\sigma_n$ the number of sites and $\beta_n$ the number of internal bonds in $\T_n$. Then
$$
\frac{\beta_n}{\sigma_n} \leq \rho_n \leq 2d-1.
$$
\label{rho}
\end{lemma}

\noindent{\it Proof.}
It suffices to prove that the configuration restricted to $\T_n$ is recurrent \cite{dhar, meester}. This can be checked in a nonambiguous way with the burning algorithm \cite{dhar}. The inequality then follows from the fact that in a recurrent, stable configuration restricted to some set, the total number of particles in the set is at least the number of internal bonds in the set, and at most $h_{max}$ times the number of sites in the set. In the case of an unstable configuration, the term ``recurrent'' is somewhat inappropriate, but we still define it to indicate configurations that pass the burning algorithm.

We use the following properties of recurrent configurations:
\begin{enumerate}
\item recurrence of a configuration restricted to a set is conserved under addition of a particle to the set, and under a legal toppling in the set.
\item if the configuration restricted to a set is recurrent, and we add an unstable site $\x$ to the set, then the configuration is recurrent restricted to the extended set.
\end{enumerate}

We will use induction in $n$.
We will show that if $\eta_n$ restricted to $\T_n$ is recurrent, then $\eta_{n+1}$  restricted to $\T_{n+1}$ is recurrent. For a starting point of the induction, we choose $n'=2d-h$, so that $|\T_{n'}|=1$, because a configuration restricted to a single site is always recurrent.

We can obtain $\eta_{n+1}$ by starting from $\eta_n$, adding a particle to the origin, and stabilizing. If we first stabilize restricted to $\T_n$, then by the first property of recurrence, the configuration restricted to $\T_n$ is still recurrent. If $\T_n = \T_{n+1}$, then we are now done. But it is also possible that at this point, there are unstable sites outside $\T_n$.

By both properties of recurrence, if we add one of these sites to $\T_n$ and stabilize with respect to the new set, the configuration restricted to the new set is still recurrent. We can repeat this step until there are no more unstable sites to be added. Then the new set is $\T_{n+1}$, and the new configuration is $\eta_{n+1}$.
\qed

\noindent
{\it Proof of Theorem \ref{topplingbounds}, part (2)}
To calculate the inequalities for $r_n$, we observe that the minimum $r_n$ would be obtained if $\V_n$ would equal the cube $\mathfrak{C}(r_n+1)$, with each site containing the maximal number of particles. The number of sites in this cube is $(2r_n+3)^d$, and the number of particles per site is then $2d-1-h$, so that
$$
n \leq (2r_n+3)^d (2d-1-h),
$$
for every $n$.

The maximal $r_n$ would be obtained if $\V_n$ would equal the diamond $\mathfrak{D}(r_n)$, containing the minimum number of particles. As we deduced that the configuration restricted to $\T_n$ is recurrent, we calculate the minimum average value of $H(\x)$ on $\mathfrak{D}(r_n-1)$, using Lemma \ref{rho}. The number of internal bonds in $\mathfrak{D}(r_n-1)$ is $(2r_n-2)^d$. The number of sites in $\mathfrak{D}(r_n-1)$ is $\frac 1d (2r_n-1)^d+o(r_n^d)$ as $r_n \to \infty$.
Thus,
$$
n \geq (2r_n-2)^d - \frac hd (2r_n-1)^d+o(r_n^d), \hspace{2cm} r_n\to \infty.
$$
As in the limit $n \to \infty$ we also have $r_n \to \infty$, this leads to the second inequality for $r_n$ in the theorem.
\qed

\medskip\noindent
The proof of Theorem \ref{topplingbounds} allows to prove an interesting characteristic of $\V_n$:

\begin{proposition}
For the sandpile model, $\V_n$ is simply connected for all $n$.
\label{simplyconnected}
\end{proposition}

\noindent{\it Proof.}
We first prove that $\T_n$ is simply connected, by contradiction.

Suppose that $\T_{n,SP,h}$ is not simply connected, that is, $\T_{n,SP,h}$ contains holes. This means there must be a site $\y$ in $\T_{n,SP,h}$ `beyond a hole', that is, a site $y$ with $T(\y)>0$, and at least one neighbor $\x$ of $\y$ with $T(\x)=0$, that is closer to the origin than $y$ itself and all neighbors $x'$ of $y$ for which $T(x')>0$.
At least one neighbor $\x'$ of $\y$ must have $T(\x')>0$, because $\T_n$ is path connected.
Thus, among the neighbors of $\y$, there must be a next-nearest neighbor pair $\z$ and $\z'$, with $d(\z') \geq d(\z)$, $T(\z)=0$ and $T(\z')>0$. But this contradicts the above derived property that for all $t$, for every next-nearest neighbor $\z$ and $\z'$ with $d(\z) \geq d(\z')$, $T^t(\z) \leq T^t(\z')$, as this property should also hold for the final configuration.

The above does not suffice to conclude that $\V_n$ is also simply connected. Since $\V_n = \T_n \cup \delta \T_n$, we must show that $T_n$ does not contain so-called {\em fjords}, i.e., places where the boundary of $\T_n$ nearly touches itself, such that $\T_n \cup \delta \T_n$ would contain a hole. But when one supposes that $\T_n$ does contain a fjord, then one can derive the same contradiction as above. Therefore, $\V_n$ is simply connected.
\qed


\begin{thebibliography}{}
\bibitem{bak} {\small\sc P. Bak, K. Tang and K. Wiesenfeld} (1988)
{Self-organized criticality},
{\em Phys. Rev. A} {\bf 38}, 364-374.

\bibitem{rossin}
Y. le Borgne and D. Rossin (2002)
On the identity of the sandpile group,
{\em Discrete Math.} {\bf 256}, 775-790.

\bibitem{ham}
A. Dartois and D. Rossin (2004)
Height Arrow Model,
{\em Formal Power Series and Algebraic Combinatorics, 16th conference, Vancouver.}

\bibitem{dhar}
{D. Dhar} (1990)
{Self-Organized Critical State of Sandpile Automaton Models},
{\em Phys. Rev. Lett.} {\bf 64}(14), 1613-1616.

\bibitem{diaconis}
{P. Diaconis and W. Fulton} (1991)
{A growth model a game, an algebra, Lagrange inversion, and characteristic classes},
{\em Rend. Sem. Mat. Univ. Pol. Torino} {\bf 49}, 95-119.

\bibitem{durrett}
{ R. Durrett} (1988)
Lecture notes on particle systems and percolation.
{\em Wadsworth and Brooks/Cole, Pacific Grove, Calif.}

\bibitem {priezzhev}
{E.V. Ivashkevich and V.B. Priezzhev}(1998)
Introduction to the sandpile model,
{\em Physica A} {\bf 254}, 97--116.

\bibitem{kleber}
{M. Kleber} (2005)
{Goldbug variations},
{\em Math. Intelligencer} {\bf 27} (1), 55-63.

\bibitem{lawler}
{G.F. Lawler, M. Bramson and D. Griffeath} (1992)
{Internal diffusion limited aggregation},
{\em Ann. Probab.} {\bf 20}(4), 2117-2140.

\bibitem{lawler2}
{G.F. Lawler} (1995)
{Subdiffusive fluctuations for internal diffusion limited aggregation},
{\em Ann. Probab.} {\bf 23}(1), 71-86.

\bibitem{levine}
{L. Levine} (2002)
{The rotor-router model},
{\em Harvard University senior thesis}.

\bibitem{peres}
{L. Levine, Y. Peres} (2005)
{Spherical Asymptotics for the Rotor-Router Model in $\mathbb{Z}^d$},
to appear in {\em Indiana University Math Journal}.

\bibitem{peres2}
{L. Levine, Y. Peres} (2007)
{Strong Spherical Asymptotics for Rotor-Router Aggregation and the Divisible Sandpile},
Preprint available at http://www.arxiv.org/abs/0704.0688.

\bibitem{liu}
S.H. Liu, T. Kaplan and L.J. Gray (1990)
Geometry and dynamics of deterministic sand piles,
{\em Phys. Rev. A} {\bf 42}(6), 3207-3212.

\bibitem{lubeck}
S. L\"{u}beck, N. Rajewsky and D.E. Wolf (2000)
A deterministic sandpile revisited,
{\em Eur. Phys. Journ. B} {\bf 13}, 715-721.

\bibitem{manna}
S.S. Manna (1991)
Two-state model of self-organized criticality,
{\em J. Phys. A: Math. Gen.} {\bf 24}, L363-369.

\bibitem{markosova}
M. Marko\v{s}ov\'{a} and P. Marko\v{s} (1992)
Analytic calculation of the attractor periods of deterministic sandpiles,
{\em Phys. Rev. A} {\bf 46}(6), 3531-3534.

\bibitem{meester}
{R. Meester, F. Redig and D. Znamenski} (2001)
{The abelian sandpile model, a mathematical introduction},
{\em Markov Proc. Rel. Fields} {\bf 7}, 509-523.

\bibitem{walkers}
{A.M. Povolotsky, V.B. Priezzhev and R.R. Shcherbakov} (1998)
{Dynamics of Eulerian walkers},
{\em Phys Rev. E} {\bf 58}, 5449-5454.

\bibitem{wiesenfeld}
K. Wiesenfeld, J. Theiler and B. McNamara(1990)
Self-organized Criticality in a Deterministic Automaton,
{\em Phys. Rev. A} {\bf 65}(8), 949-952.
\end{thebibliography}
\end{document}